\documentclass[MSNbibl,number,citesort,seceqn,dvips]{arxbj}
\usepackage{dcolumn,upgreek}
\usepackage{graphicx}

\aid{0}
\volume{19}
\issue{1}
\pubyear{2013}
\firstpage{205}
\lastpage{227}
\doi{10.3150/11-BEJ399}

\makeatletter
\newcolumntype{d}[1]{D{.}{.}{#1}}
\newtheorem{theorem}{Theorem}[section]
\newremark{rem}{Remark}[section]
\newproclaim{assumption}{Assumption}[section]
\newtheorem{proposition}{Proposition}[section]
\newremark{example}{Example}[section]

\newcommand{\field}[1]{\mathbb{#1}}
\newcommand{\R}{\field{R}}
\newcommand{\p}{\field{P}}
\newcommand{\Z}{\field{Z}}
\def\argmax{\operatorname{argmax}}
\makeatother

\begin{document}
\begin{frontmatter}

\title{Inference for modulated stationary processes}
\runtitle{Modulated stationary processes}

\begin{aug}
\author{\fnms{Zhibiao} \snm{Zhao}\thanksref{e1}\corref{}\ead[label=e1,mark]{zuz13@stat.psu.edu}}%
\and
\author{\fnms{Xiaoye} \snm{Li}\thanksref{e2}\ead[label=e2,mark]{xul117@stat.psu.edu}}
\runauthor{Z. Zhao and X. Li} 
\address{Department of Statistics, Penn State University, University
Park, PA
16802, USA.\\
\printead{e1,e2}}
\end{aug}

\received{\smonth{1} \syear{2011}}
\revised{\smonth{8} \syear{2011}}

%
\begin{abstract}
We study statistical inferences for a class of modulated stationary processes
with time-dependent variances. Due to
non-stationarity and the large number of unknown parameters,
existing methods for stationary, or locally
stationary, time series are not applicable. Based on a
self-normalization technique, we address several inference
problems, including a self-normalized central limit theorem,
a self-normalized cumulative sum test for the change-point problem,
a long-run variance estimation through blockwise
self-normalization, and a self-normalization-based wild bootstrap.
Monte Carlo simulation studies
show that the proposed self-normalization-based methods
outperform stationarity-based alternatives. We demonstrate the
proposed methodology using two real data sets: annual mean
precipitation rates in Seoul from 1771--2000, and quarterly
U.S. Gross National Product growth rates from 1947--2002.
\end{abstract}

%

\end{frontmatter}

\section{Introduction}\label{sec:intr}

In time series analysis, stationarity requires that dependence
structure be sustained over time, and thus we can borrow
information from one time period to study model dynamics over
another period; see Fan and Yao~\cite{fan} for nonparametric
treatments and Lahiri~\cite{lahiri2} for various resampling and block
bootstrap methods. In practice, however, many climatic,
economic and financial time series are non-stationary and therefore
challenging to analyze. First, since dependence structure
varies over time, information is more localized. Second,
non-stationary processes often require extra parameters to
account for time-varying structure. One way to overcome these issues
is to impose certain local stationarity; see, for
example, Dahlhaus~\cite{dahlhaus} and Adak~\cite{adak} for
spectral representation frameworks and Dahlhaus and Polonik~\cite{dahlhaus2} for a time domain
approach.

In this article we study a class of modulated stationary processes (see
Adak~\cite{adak})
%
\begin{equation}\label{eq:xinons}
X_i = \mu+ \sigma_i e_i,\qquad  i=1,\ldots,n,
\end{equation}
where $e_i$ are stationary time series with zero mean, and
$\sigma_i>0$ are unknown constants adjusting for time-dependent
variances. Then $X_i$ oscillates around the constant mean $\mu$,
whereas its variance changes over time in an unknown manner. In
the special case of $\sigma_i\equiv1$, (\ref{eq:xinons})
reduces to stationary case. If $\sigma_i=s(i/n)$ for a
Lipschitz continuous function $s(t)$ on $[0,1]$, then
(\ref{eq:xinons}) is locally stationary.
For the general
non-stationary case (\ref{eq:xinons}), the number of unknown
parameters is larger than the number of observations, and
it is infeasible to estimate $\sigma_i$. Due to
non-stationarity and the large number of unknown
parameters, existing methods that are developed for (locally)
stationary processes are not applicable, and our main purpose
is to develop new statistical inference techniques.\looseness=1

First, we establish a uniform strong approximation result which
can be used to derive a self-normalized central limit theorem
(CLT) for the sample mean $\bar{X}$ of (\ref{eq:xinons}). For
stationary case $\sigma_i\equiv1$, by Fan and Yao~\cite{fan},
under mild mixing conditions,
%
\begin{equation}\label{eq:fanyao}
\sqrt{n}(\bar{X}-\mu)\Rightarrow N(0,\tau^2),
 \qquad\mbox{where }  \tau^2=\gamma_0+2\sum^\infty_{k=1}\gamma_k
 \mbox{ and }  \gamma_k =\operatorname{Cov}(e_i,e_{i+k}).
\end{equation}
For
the modulated stationary case (\ref{eq:xinons}), it is non-trivial
whether $\sqrt{n}(\bar{X}-\mu)$ has a CLT without imposing
further assumptions on $\sigma_i$ and the dependence structure
of $e_i$. Moreover, even when the latter CLT exists, it is
difficult to estimate the limiting variance due to the large
number of unknown parameters; see De Jong and Davidson~\cite{dejong} for related work
assuming a near-epoch dependent mixing framework.
Zhao~\cite{zhao} studied confidence interval construction for $\mu$ in~(\ref
{eq:xinons})
by assuming a block-wise asymptotically equal cumulative variance
assumption. The latter assumption is rather restrictive and
essentially requires that block averages be asymptotically independent
and identically distributed (i.i.d.). In this article, we
deal with the more general setting~(\ref{eq:xinons}).
Under a strong invariance
principle assumption, we establish a self-normalized CLT with
the self-normalizing constant adjusting for time-dependent
non-stationarity. The obtained CLT is an extension of the
classical CLT for i.i.d.
data or stationary time series to modulated stationary processes.
Furthermore, we extend the idea to linear combinations of means
over different time periods, which allows us to address
inference regarding mean levels over multiple time periods.

Second, we study the wild bootstrap for modulated stationary processes.
Since the seminal work of Efron~\cite{efron}, a great deal of research has
been done
on the bootstrap under various settings, ranging from
bootstrapping for i.i.d. data in Efron~\cite{efron}, wild bootstrapping for
independent observations with possibly non-constant variances
in Wu~\cite{wuj} and Liu~\cite{liu}, to various block bootstrapping and
resampling methods for stationary time series in K\"{u}nsch~\cite{kunsch},
Politis and Romano~\cite{politis}, B\"{u}hlmann~\cite{buhlmann} and the
monograph Lahiri~\cite{lahiri2}. With the established self-normalized
CLT, we propose a wild bootstrap procedure that is tailored to
deal with modulated stationary processes: the dependence is
removed through a scaling factor, and the non-constant variance
structure of the original data is preserved in the wild
bootstrap data-generating mechanism.
Our simulation study shows that the wild bootstrap
method outperforms the widely used stationarity-based block bootstrap.

Third, we address change-point analysis. The change-point problem
has been an active area of research; see Pettitt~\cite{pettitt} for
proportion changes in binary data, Horv\'ath~\cite{horvath} for mean
and variance changes in Gaussian observations,
Bai and Perron~\cite{bai} for coefficient changes in linear models, Aue \textit{et al.}~\cite{aue2008a}
for coefficient changes in polynomial regression with
uncorrelated errors, Aue \textit{et al.}~\cite{aue2008b} for mean change in time series
with stationary errors,
Shao and Zhang~\cite{shaox} for change-points for stationary time series and the
monograph by Cs\"org\H{o} and Horv\'ath~\cite{csorgo} for more
discussion. Most of these works deal with
stationary and/or independent data. Hansen~\cite{hansen}
studied tests for constancy of parameters in linear regression
models with non-stationary regressors and conditionally
homoscedastic martingale difference errors. Here we consider
%
\begin{equation}\label{eq:null}
H_0\dvt   X_i=\mu_i+\sigma_ie_i, \mu_1 = \cdots= \mu_n,\qquad
  H_a\dvt   \mu_1 = \cdots= \mu_{J} \not=
\mu_{J+1} = \cdots= \mu_n,
\end{equation}
where $J$ is an unknown change point. The aforementioned works mainly
focused on detecting changes in mean
while the error variance is constant. On the other hand, researchers
have also realized the importance of the variance/covariance structure in
change point analysis. For example, Incl\'an and Tiao~\cite{inclan} studied change in
variance for independent data,
and Aue \textit{et al.}~\cite{aue2009} and Berkes, Gombay and Horv\'ath~\cite{berkes} considered change in covariance
for time series data.
To our knowledge, there
has been almost no attempt to advance change point analysis under the
non-constant variances framework in (\ref{eq:null}).
Andrews~\cite{andrews} studied change point problem under near-epoch dependence
structure that allows for non-stationary processes, but his
Assumption~1(c) on page 830 therein essentially implies that the
process has constant variance.
The popular cumulative
sum (CUSUM) test is developed for stationary time series and
does not take into account the time-dependent variances. Using
the self-normalization idea, we propose a self-normalized CUSUM
test and a wild bootstrap method to obtain its critical value.
Our empirical studies show that the usual CUSUM tests tend to
over-reject the null hypothesis in the presence of non-constant
variances. By contrast, the self-normalized CUSUM test yields
size close to the nominal level.

Fourth, we estimate the long-run variance $\tau^2$ in
(\ref{eq:fanyao}). Long-run variance plays an essential role
in statistical inferences involving time series. Most works
in the literature deal with stationary processes through various block bootstrap
and subsampling approaches; see Carlstein~\cite{carlstein}, K\"{u}nsch~\cite{kunsch},
Politis and Romano~\cite{politis}, G\"{o}tze and K\"{u}nsch~\cite{gotze} and the monograph
Lahiri~\cite{lahiri2}. De Jong and Davidson~\cite{dejong} established the consistency of kernel
estimators of
covariance matrices under a near epoch dependent mixing condition.
Recently, M\"uller~\cite{muller} studied robust long-run
variance estimation for locally stationary process. For
model (\ref{eq:xinons}), the error process
$\{e_i\}$ is contaminated with unknown standard deviations
$\{\sigma_i\}$, and we apply blockwise
self-normalization to remove non-stationarity, resulting in
asymptotically stationary blocks.

Fifth, the proposed methods can be
extended to deal with the linear regression model
%
\begin{equation}\label{eq:lr}
X_i = U_i \beta+ \sigma_i e_i,
\end{equation}
where $U_i=(u_{i,1},\ldots,u_{i,p})$ are deterministic
covariates, and $\beta=(\beta_1,\ldots,\beta_p)'$ is the unknown
column vector of parameters. For $p=2$,
Hansen~\cite{hansen1995} established the asymptotic normality of the
least-squares estimate
of the slope parameter under a fairly general framework of
non-stationary errors.
While Hansen~\cite{hansen1995} assumed that the errors form a martingale
difference array so that they are uncorrelated,
the framework in (\ref{eq:lr}) is more general in that it allows for
correlations. On the other hand,
Hansen~\cite{hansen1995} allowed the conditional volatilities to
follow an autoregressive model, hence introducing stochastic volatilities.
Phillips, Sun and Jin~\cite{phi} considered (\ref{eq:lr}) for stationary
errors, and their approach is not applicable here due to the
unknown non-constant variances $\sigma^2_i$. In Section~\ref{sec:ext}
we consider self-normalized CLT for
the least-squares estimator of $\beta$ in (\ref{eq:lr}). In the
polynomial regression
case $u_{i,r}=(i/n)^{r-1}$, Aue \textit{et al.}~\cite{aue2008a} studied a likelihood-based
test for constancy
of $\beta$ in (\ref{eq:lr}) for uncorrelated errors with constant variance.
Due to the presence of correlation and time-varying variances, it is
more challenging to study the change point problem
for (\ref{eq:lr}) and this is beyond the scope of this
article.\looseness=1

The rest of this article is organized as follows. We present
theoretical results in Section~\ref{sec:main}. Sections
\ref{sec:simu}--\ref{sec:app} contain Monte Carlo studies and
applications to two real data sets.

\section{Main results}\label{sec:main}
For sequences $\{a_n\}$ and $\{b_n\}$, write $a_n=\mathrm{O}(b_n)$,
$a_n=\mathrm{o}(b_n)$ and $a_n\asymp b_n$, respectively, if
$|a_n/b_n|<c_1$, $a_n/b_n\to0$ and $c_2 < |a_n/b_n| <c_3$,
for some constants $0<c_1,c_2,c_3<\infty$. For $q>0$ and a random
variable $e$, write $e\in\mathcal{L}^q$
if $\|e\|_q:=\{E(|e|^q)\}^{1/q}<\infty$.

\subsection{Uniform approximations for modulated stationary
processes}\label{sec:con}

In (\ref{eq:xinons}), assume without loss of generality that
$E(e_i)=0$ and $E(e^2_i)=1$ so that $\{e_i\}$ and $\{e^2_i-1\}$
are centered stationary processes. With the convention
$S_0=S^*_0=0$, define
%
\begin{equation}\label{eq:SnSn}
S_i = \sum^i_{j=1} e_j  \quad \mbox{and}\quad
  S^*_i = \sum^i_{j=1} (e^2_j-1), \qquad   i=1,2,\ldots.
\end{equation}

\begin{assumption}\label{assump:1}
There exist standard Brownian motions $\{B_t\}$ and
$\{B^*_t\}$ such that
%
\begin{equation}\label{eq:sip}
\max_{1\le i\le n}|S_i - \tau B_i |=\mathrm{o}_\mathrm{a.s.}(\Delta_n)\quad  \mbox
{and}\quad
\max_{1\le i\le n}|S^*_i - \tau^* B^*_i |=\mathrm{o}_\mathrm{a.s.}(\Delta_n),
\end{equation}
where $\Delta_n$ is the approximation error, $\tau^2$ and
$\tau^{*2}$ are the long-run variances of $\{e_i\}$ and
$\{e^2_i-1\}$, respectively. Further assume $\tau^2>0$ to avoid
the degenerate case $\tau^2=0$.
\end{assumption}

The uniform approximations in (\ref{eq:sip}) are generally
called strong invariance principle. The two Brownian motions $\{B_t\}$ and
$\{B^*_t\}$ may be defined on different probability spaces and hence
are not jointly distributed,
which is not an issue because our argument does not depend on their
joint distribution.
To see how to use (\ref{eq:sip}), under $H_0$ in
(\ref{eq:null}), consider
%
\begin{equation}\label{eq:Fj}
F_j=j (\underline{X}_j-\mu)  \quad\mbox{and}\quad
\underline{V}^2_j=\sum^j_{i=1}(X_i-\underline{X}_j)^2,\qquad \mbox
{where }
\underline{X}_j=j^{-1}\sum^j_{i=1}X_i.
\end{equation}
Theorem~\ref{thm:0} below presents uniform approximations for
$F_j$ and $V^2_j$. Define
%
\begin{eqnarray}\label{eq:volwei}
r_n &= &|\sigma_n| + \sum^n_{i=2} |\sigma_i-\sigma_{i-1}|
 \quad \mbox{and}\quad
r^*_n = |\sigma^2_n| + \sum^n_{i=2} |\sigma^2_i-\sigma^2_{i-1}|,
\\[-3pt]
\label{eq:Omegan}
\Sigma^2_j&=&\sum^j_{i=1} \sigma^2_i   \quad\mbox{and}\quad
\Sigma^{*2}_j=\Biggl( \sum^j_{i=1}\sigma^4_i \Biggr)^{1/2}.\vspace*{-2pt}
\end{eqnarray}

\begin{theorem}\label{thm:0}
Let (\ref{eq:sip}) hold. For any $c\in(0,1]$, the following
uniform approximations hold:
%
\begin{eqnarray}
&&\max_{cn \le j\le n} \Biggl|F_j-\tau\sum^j_{i=1}
\sigma_i (B_i-B_{i-1}) \Biggr|=\mathrm{O}_\mathrm{ a.s.}(r_n\Delta_n),\label{eq:thm0a}
\\[-3pt]
&&\max_{cn \le j\le n} |\underline{V}^2_j - \Sigma^2_j |=
\mathrm{O}_\mathrm{ p}\{(r^2_n\Delta^2_n + \Sigma^2_n)/n + \Sigma^{*2}_n+r^*_n \Delta
_n\}.\label{eq:thm0b}\vspace*{-2pt}
\end{eqnarray}
\end{theorem}

Theorem~\ref{thm:0} provides quite general results under (\ref{eq:sip}).
We now discuss sufficient conditions for (\ref{eq:sip}). Shao~\cite{shao}
obtained sufficient mixing conditions for
(\ref{eq:sip}). In this article, we briefly introduce the
framework in Wu~\cite{wuw}. Assume that $e_i$ has the causal
representation
$e_i = G(\ldots,\varepsilon_{i-1}, \varepsilon_i)$,
where
$\varepsilon_i$ are i.i.d. innovations,
and $G$ is a measurable function such that $e_i$ is
well defined. Let $\{\varepsilon_i'\}_{i\in\Z}$ be an independent copy of
$\{\varepsilon_i\}_{i\in\Z}$.
Assume
%
\begin{equation}\label{eq:pro1con}
\sum_{i=1}^\infty i \|e_i - e'_i\|_8< \infty, \qquad\mbox{where }
  e'_i = G(\ldots, \varepsilon_{-1}, \varepsilon_0',
\varepsilon_1, \ldots, \varepsilon_{i-1}, \varepsilon_i).
\end{equation}
Proposition~\ref{pro:1} below follows from
Corollary 4 in Wu~\cite{wuw}.\vspace*{-2pt}

\begin{proposition}\label{pro:1}
Assume that (\ref{eq:pro1con}) holds.
Then (\ref{eq:sip}) holds with $\Delta_n=n^{1/4}\log(n)$, the
optimal rate up to a logarithm factor.\vspace*{-2pt}
\end{proposition}

For linear process $e_i = \sum^\infty_{j=0} a_j
\varepsilon_{i-j}$ with $\varepsilon_i\in\mathcal{L}^8$ and
$E(\varepsilon_i)=0$,
$\|e_i-e'_i\|_8=\|\varepsilon_0-\varepsilon'_0\|_8 |a_i|$. If
$\sum^\infty_{i=1} i|a_i|<\infty$, then (\ref{eq:sip}) holds
with $\Delta_n=n^{1/4}\log(n)$. For many nonlinear time series,
$\|e_i - e'_i\|_8$ decays exponentially fast and hence
(\ref{eq:pro1con}) holds; see
Section 3.1 of Wu~\cite{wuw}. From now on we assume (\ref{eq:sip})
holds with $\Delta_n=n^{1/4}\log(n)$.\vspace*{-2pt}

\begin{rem}\label{rmk:moment}
If $e_i$ are i.i.d. with $E(e_i)=0$ and $e_i\in\mathcal{L}^q$ for some
$2<q\le4$,
the celebrated ``Hungarian embedding'' asserts
that $\sum^i_{j=1}e_j$ satisfies a strong invariance principle with the
optimal rate $\mathrm{o}_\mathrm{ a.s.}(n^{1/q})$. Thus, it is
necessary to have
the moment assumption $e_i\in\mathcal{L}^8$ in Proposition~\ref{pro:1} in
order to ensure
strong invariance principles for both $S_i$ and $S^*_i$ in (\ref
{eq:SnSn}) with
approximation rate $n^{1/4}\log(n)$. On the other hand, one can relax
the moment assumption by loosening the
approximation rate. For example, by Corollary 4 in Wu~\cite{wuw}, assume
$e_i\in\mathcal{L}^{2q}$ for some $q>2$ and $\sum_{i=1}^\infty i \|e_i -
e_i^*\|_{2q}< \infty$, then (\ref{eq:sip}) holds with
$\Delta_n=n^{1/\min(q,4)}\log(n)$.\vspace*{-2pt}
\end{rem}

As shown in Examples~\ref{exmp:3}--\ref{exmp:4} below, $r_n$
and $r^*_n$ in (\ref{eq:volwei}) often have tractable bounds.\vadjust{\goodbreak}

\begin{example}\label{exmp:3}
If $\sigma_i$ is non-decreasing in $i$, then $\sigma_n\le
r_n\le2\sigma_n$ and $\sigma^2_n\le r^*_n\le2\sigma^2_n$. If
$\sigma_i$ is non-increasing in $i$, then $r_n = \sigma_1$ and
$r^*_n=\sigma^2_1$. If $\sigma_i$ are piecewise constants with finitely
many pieces,
then $r_n,r^*_n=\mathrm{O}(1)$.
\end{example}

\begin{example}\label{exmp:6}
Let $\sigma_i=s(i/n^\gamma)$ for $\gamma\in[0,1]$ and a
Lipschitz continuous function $s(t), t\in[0,\infty),
\sup_{t\in[0,\infty)}s(t)<\infty$. Then
$r_n,r^*_n=\mathrm{O}(n^{1-\gamma})$. If $\gamma=1$, we obtain a locally
stationary case with the time window $i/n\in[0,1]$; if $\gamma\in
[0,1)$, we have the infinite time window $[0,\infty)$ as
$n/n^\gamma\to\infty$, which may be more reasonable for data
with a long time horizon.
\end{example}

\begin{example}\label{exmp:4}
If $\sigma_i=i^\beta L(i)$ for a slowly varying function
$L(\cdot)$ such that $L(cx)/L(x)\to1$ as $x\to\infty$ for all
$c>0$. Then we can show $r_n=\mathrm{O}\{n^\beta L(n)\}$ or $\mathrm{O}(1)$ and $r^*_n=\mathrm{O}\{
n^{2\beta} L^2(n)\}$
or $\mathrm{O}(1)$, depending on whether
$\beta>0$ or $\beta<0$. For the boundary case $\beta=0$, assume
$L(i+1)/L(i)=1+\mathrm{O}(1/i)$ uniformly, then
$r_n=L(n)+\mathrm{O}(1) \sum^n_{i=2} L(i)/i=\mathrm{O}\{\log(n) \times\max_{1\le i\le n} L(i)\}
$. Similarly, $r^*_n=\mathrm{O}\{\log(n)  \max_{1\le i\le n} L^2(i)\}$.
\end{example}

\subsection{Self-normalized central limit
theorem}\label{sec:cltx}

In this section we establish a self-normalized
CLT for the sample average $\bar{X}$. To understand how
non-stationarity makes
this problem difficult, elementary calculation shows
%
\begin{equation}\label{eq:lrvnon}
\operatorname{Var}\bigl\{\sqrt{n}(\bar{X}-\mu)\bigr\} = \frac{\gamma_0}{n} \sum^n_{i=1}
\sigma^2_i +
\frac{2}{n} \sum_{1\le i<j\le n} \sigma_i\sigma_j \gamma_{j-i}
:=\tau^2_n,
\end{equation}
where $\gamma_k=\operatorname{ Cov}(e_0,e_k)$. In the stationary case $\sigma
_i\equiv1$, under condition
$\sum^\infty_{k=0}|\gamma_k|<\infty$, $\tau^2_n\to\tau^2$, the
long-run variance in (\ref{eq:fanyao}). For non-constant variances,
it is difficult to deal with $\tau^2_n$ directly, due to
the large number of unknown parameters and complicated
structure. See De Jong and Davidson~\cite{dejong} for a kernel estimator of $\tau^2_n$ under
a near-epoch dependent mixing framework.

To attenuate the aforementioned issue, we apply the uniform
approximations in Theorem~\ref{thm:0}. Assume that
(\ref{eq:thmtestcon}) below holds. Note that the increments
$B_i-B_{i-1}$ of standard Brownian motions are i.i.d. standard normal
random variables. By (\ref{eq:thm0a}),
$n(\bar{X}-\mu)$ is equivalent to $N(0,\tau^2\Sigma^2_n)$ in
distribution. By (\ref{eq:thm0b}), $\underline{V}_n/\Sigma_n\to
1$ in probability. By Slutsky's theorem, we have Proposition~\ref{cor:1}.

\begin{proposition}\label{cor:1}
Let (\ref{eq:sip}) hold with $\Delta_n= n^{1/4}\log(n)$. For
$r_n, r^*_n, \Sigma^2_n, \Sigma^{*2}_n$ in (\ref{eq:volwei})--(\ref{eq:Omegan}), assume
%
\begin{equation}\label{eq:thmtestcon}
\delta_n=r_n\Delta_n/\Sigma_n +(r^*_n \Delta_n +
\Sigma^{*2}_n)/\Sigma^2_n\to0.
\end{equation}
Recall $\underline{V}_n^2$ in (\ref{eq:Fj}). Then as
$n\to\infty$, $n(\bar{X}-\mu)/\underline{V}_n\Rightarrow
N(0,\tau^2)$.
Consequently, a $(1-\alpha)$ asymptotic confidence interval for
$\mu$ is $\bar{X} \pm z_{\alpha/2} \hat\tau\underline{V}_n/n$,
where $\hat\tau$ is a consistent estimate of $\tau$ (Section~\ref{sec:lrv} below), and $z_{\alpha/2}$ is $(1-\alpha/2)$
standard normal quantile.
\end{proposition}

Proposition~\ref{cor:1} is an extension of the classical CLT
for i.i.d. data or stationary processes to modulated stationary processes.
If $X_i$ are i.i.d., then
$n(\bar{X}-\mu)/\underline{V}_n\Rightarrow N(0,1)$. In
Proposition~\ref{cor:1}, $\tau^2$ can be viewed as the variance
inflation factor due to the dependence of $\{e_i\}$. For
stationary data, the sample variance $\underline{V}^2_n/n$ is a
consistent estimate of the population variance. Here, for
non-constant variances case (\ref{eq:xinons}), by (\ref{eq:thm0b}) in
Theorem~\ref{thm:0}, $\underline{V}^2_n/n$ can be viewed as an
estimate of the time-average ``population variance''
$\Sigma^2_n/n$. So, we can interpret the CLT in Proposition~\ref{cor:1} as a self-normalized CLT for modulated stationary processes
with the self-normalizing term $\underline{V}_n$,
adjusting for non-stationarity due to
$\sigma_1,\ldots,\sigma_n$ and $\tau^2$, accounting for
dependence of $\{e_i\}$. Clearly, parameters
$\sigma_1,\ldots,\sigma_n$ are canceled out through
self-normalization. Finally, condition (\ref{eq:thmtestcon}) is
satisfied in
Example~\ref{exmp:6} with $\gamma>3/4$ and Example~\ref{exmp:4}
with $\beta>-1/4$.

In classical statistics, the width of confidence intervals usually
shrinks as sample size increases. By Proposition~\ref{cor:1}
and Theorem~\ref{thm:0}, the width of the constructed
confidence interval for $\mu$ is proportional to
$\underline{V}_n/n$ or, equivalently, $\Sigma_n/n$. Thus, a
necessary and sufficient condition for shrinking confidence
interval is $\sum^n_{i=1}\sigma^2_i/n^2\to0$, which is
satisfied if $\sigma_i=\mathrm{o}(\sqrt{i})$. An intuitive explanation
is as follows. For i.i.d. data, sample mean converges at a rate of
$\mathrm{O}(\sqrt{n})$. In (\ref{eq:xinons}), if $\sigma_i$ grows
faster than $\mathrm{O}(\sqrt{i})$, the contribution of a new
observation is negligible relative to its noise level.

\begin{example}\label{exmp:ci}
If $\sigma_i\asymp i^{\beta}$ with $\beta\in[0,1/2)$, the
length of confidence interval is proportional to
$\Sigma_n/n\asymp n^{\beta-1/2}$. In particular, if
$c_1<\sigma_i<c_2$ for some positive constants $c_1$ and $c_2$,
then $\Sigma_n/n$ achieves the optimal rate $\mathrm{O}(n^{-1/2})$. If
$\sigma_i\asymp\log(i)$, then $\Sigma_n/n\asymp
\log(n)/\sqrt{n}$.
\end{example}

The same idea can be extended to linear combinations of means
over multiple time periods. Suppose we have observations
from $k$ consecutive time periods $\mathcal{T}_1,\ldots,\mathcal{T}_k$, each of the form~(\ref{eq:xinons}) with different means,
denoted by $\mu_1,\ldots,\mu_k$, and each having time-dependent variances.
Let $\nu=\beta_1\mu_1+\cdots+\beta_k\mu_k$ for given
coefficients $\beta_1,\ldots,\beta_k$. For example, if we are
interested in mean change from $\mathcal{T}_1$ to $\mathcal{T}_2$, we
can take $\nu=\mu_2-\mu_1$; if we are interested in whether the
increase from $\mathcal{T}_3$ to $\mathcal{T}_4$ is larger than that
from $\mathcal{T}_1$ to $\mathcal{T}_2$, we can let
$\nu=(\mu_4-\mu_3)-(\mu_2-\mu_1)$. Proposition~\ref{thm:lcm}
below extends Proposition~\ref{cor:1} to multiple means.

\begin{proposition}\label{thm:lcm}
Let $\nu=\beta_1\mu_1+\cdots+\beta_k\mu_k$. For $\mathcal{T}_j$,
denote its sample size $n_j$ and its sample
average $\bar{X}(j)$. Assume that (\ref{eq:thmtestcon}) holds for each
individual time period $\mathcal{T}_j$ and, for simplicity, that
$n_1,\ldots,n_k$ are of the same order. Then
\[
\frac{\sum^k_{j=1} \beta_j \bar{X}(j)-\nu}{ \Lambda_n}\Rightarrow
N(0,\tau^2),
 \qquad\mbox{where }  \Lambda^2_n = \sum^k_{j=1}
\biggl\{ \frac{\beta^2_j}{n^2_j } \sum_{i\in\mathcal{T}_j} [X_i-\bar
{X}(j)]^2 \biggr\}.
\]
\end{proposition}

\subsection{Wild bootstrap for self-normalized
statistic}\label{sec:wild} Recall $\sigma_i e_i$ in
(\ref{eq:xinons}). Suppose we are interested in the
self-normalized statistic
\[
H_n=\frac{\sum^n_{i=1}\sigma_i e_i}{\sqrt{\sum^n_{i=1} \sigma^2_i e^2_i}}.\vadjust{\goodbreak}
\]
For problems with small sample sizes, it is natural to use
bootstrap distribution instead of the convergence
$H_n\Rightarrow N(0,\tau^2)$ in Proposition~\ref{cor:1}. Wu~\cite{wuj} and
Liu~\cite{liu} have pioneered the work on the wild bootstrap
for independent data with non-identical distributions. We shall
extend their wild bootstrap procedure to the modulated stationary
process (\ref{eq:xinons}).

Let $\{\alpha_i\}$ be i.i.d. random variables independent of
$\{e_i\}$ satisfying $\alpha_i\in\mathcal{L}^3, E(\alpha_i)=0,
E(\alpha^2_i)=1$. Define the self-normalized statistic based on the
following new
data:
\[
H^*_n=\frac{\sum^n_{i=1}\xi_i}{\sqrt{\sum^n_{i=1} (\xi_i-\bar{\xi})^2}},
\qquad \mbox{where }  \xi_i = \sigma_ie_i\alpha_i  \mbox{ and }
\bar{\xi}=\frac{\xi_1+\cdots+\xi_n}{n}.
\]
Clearly, $\xi_i$ inherits the non-stationarity structure of
$\sigma_i e_i$ by writing $\xi_i=\sigma_i e^*_i$ with
$e^*_i=e_i\alpha_i$. On the other hand, for the new error
process $\{e^*_i\}$, $E(e^{*2}_i)=E(e^2_i)=1$ and $\operatorname{
Cov}(e^*_i,e^*_j)=0$ for $i\ne j$. Thus, $\{e^*_i\}$ is a white
noise sequence with long-run variance one. By Proposition
\ref{cor:1}, the scaled version $H_n/\tau\Rightarrow N(0,1)$ is
robust against the dependence structure of $\{e_i\}$, so
we expect that $H_n^*$ should be close to
$H_n/\tau$ in distribution.

\begin{theorem}\label{thm:bootstrap}
Let the conditions in Proposition~\ref{cor:1} hold. Further
assume
%
\begin{equation}\label{eq:bootcon}
\Biggl(\sum^n_{i=1} \sigma^3_i \Biggr)^2 \Biggl(\sum^n_{i=1} \sigma^2_i
\Biggr)^{-3} \to0.
\end{equation}
Let $\hat\tau$ be a consistent estimate of $\tau$. Denote by
$\p^*$ the conditional law given $\{e_i\}$. Then
%
\begin{equation}\label{eq:bootstrap}
\sup_{x\in\R} |\p^*(H^*_n\le x) - \p(H_n/\hat\tau\le x) | \to0,
\qquad\mbox{in probability}.
\end{equation}
\end{theorem}

Theorem~\ref{thm:bootstrap} asserts that, $H^*_n$ behaves like the scaled
version $H_n/\hat\tau$, with the scaling factor $\hat\tau$
coming from the dependence of $\{e_i\}$. Here we use the sample
mean $\bar{X}$ in (\ref{eq:xinons}) to illustrate a wild
bootstrap procedure to obtain the distribution of
$n(\bar{X}-\mu)/(\tau\underline{V}_n)$ in Proposition~\ref{cor:1}.
\begin{enumerate}[(iii)]
\item[(i)] Apply the method in Section~\ref{sec:lrv} to
$X_1,\ldots,X_n$ to obtain a consistent estimate
$\hat\tau$ of $\tau$.
\item[(ii)] Subtract the sample mean $\bar{X}$ from data to
obtain $\epsilon_i=X_i-\bar{X}, i=1,\ldots,n$.

\item[(iii)] Generate i.i.d. random variables
$\alpha_1,\ldots,\alpha_n$ satisfying $E(\alpha_i)=0,
E(\alpha^2_i)=1$.

\item[(iv)] Based on $\epsilon_i$ in (ii) and $\alpha_i$ in
(iii), generate bootstrap data $\xi^b_i = \epsilon_i
\alpha_i$, and compute
\begin{eqnarray*}
H^b_n=\frac{\sum^n_{i=1}\xi^b_i}{\hat\tau^b \sqrt{\sum^n_{i=1} (\xi
^b_i-\bar{\xi}^b)^2}},
\end{eqnarray*}
where $\hat\tau^b$ is a long-run variance estimate (see
Section~\ref{sec:lrv}) for bootstrap data $\xi^b_i$.

\item[(v)] Repeat (iii)--(iv) many times and
use the empirical distribution of those realizations of
$H^b_n$ as the distribution of $n(\bar{X}-\mu)/(\tau
\underline{V}_n)$.
\end{enumerate}

The proposed wild bootstrap is an extension of that in Liu~\cite{liu} for
independent data to
modulated stationary case, and
it has two appealing features. First, the scaling factor
$\hat\tau$ makes the statistic independent of the dependence
structure. Second, the bootstrap data-generating mechanism is
adaptive to unknown time-dependent variances $\{\sigma^2_i\}$.
For the distribution of $\alpha_i$ in step~(iii), we use
$\p(\alpha_i=-1)=\p(\alpha_i=1)=1/2$, which has some desirable
properties. For example, it preserves the magnitude and range of the
data. As shown by Davidson and Flachaire~\cite{davidson}, for certain hypothesis testing
problems in linear regression
models with symmetrically distributed errors, the bootstrap
distribution is exactly equal
to that of the test statistic; see Theorem 1 therein.

For the purpose of comparison, we briefly introduce the widely
used block bootstrap for a stationary time series $\{X_i\}$
with mean $\mu$. By (\ref{eq:fanyao}),
$\sqrt{n}(\bar{X}-\mu)\Rightarrow N(0,\tau^2)$. Suppose that we
want to bootstrap the distribution of $\sqrt{n}(\bar{X}-\mu)$.
Let $k_n, \ell_n, \mathcal{I}_1,\ldots,\mathcal{I}_{\ell_n}$ be
defined as in Section~\ref{sec:lrv} below. The non-overlapping
block bootstrap works in the following way:
\begin{enumerate}[(iii)]
\item[(i)] Take a simple random sample of size $\ell_n$
with replacement from the blocks $\mathcal{I}_1,\ldots,\mathcal{I}_{\ell_n}$, and form the bootstrap
data $X^b_1,\ldots,X^b_{n'}, n'=k_n\ell_n,$ by pooling
together $X_i$s for which the index $i$ is within
those selected blocks.
\item[(ii)] Let $\bar{X}^b$ be the sample average of
$X^b_1,\ldots,X^b_{n'}$. Compute $\Xi_n = \sqrt{n'} \{
\bar{X}^b - E^*( \bar{X}^b ) \}$, where $E^*(\bar{X}^b)
= \sum^{n'}_{i=1}X_i/n'$ is the conditional expectation
of $\bar{X}^b$ given $\{X_i\}$.
\item[(iii)] Repeat (i)--(ii) many times and
use the empirical distribution of $\Xi_n$'s as the
distribution of $\sqrt{n}(\bar{X}-\mu)$.
\end{enumerate}

In step (ii), another choice is the studentized
version $\tilde\Xi_n=\sqrt{n'} \{
\bar{X}^b - E^*( \bar{X}^b ) \}/\hat\tau^b$, where
$\hat\tau^b$ is a consistent estimate of $\tau$ based on
bootstrap data. Assuming stationarity and $k_n\to\infty$,
the blocks are asymptotically independent and share the same model
dynamics as the whole data, which validates the above block bootstrap. We
refer the reader to Lahiri~\cite{lahiri2} for detailed discussions.
For a non-stationary process, block bootstrap
is no longer valid, because individual blocks are not
representative of the whole data. By contrast, the proposed wild
bootstrap is adaptive to unknown dependence and
the non-constant variances structure.

\subsection{Change point analysis: Self-normalized CUSUM test}\label
{sec:cusum} To test a
change point in the mean of a process $\{X_i\}$, two popular
CUSUM-type tests (see Section~3 of Robbins \textit{et al.}~\cite{robbins} for a review and
related references) are
%
\begin{equation}\label{eq:cusum}
T^1_n = \max_{c n \le j\le(1-c) n}\frac{\hat\tau^{-1} |S_X(j)|}{\sqrt
{j(1-j/n)}}
\quad  \mbox{and}\quad
T^2_n = \max_{c n \le j\le(1-c) n}\hat\tau^{-1} |S_X(j)|,
\end{equation}
where $\hat\tau^2$ is a consistent estimate of the long-run
variance $\tau^2$ of $\{X_i\}$, and
%
\begin{equation}\label{eq:SXj}
S_X(j)
= \biggl( 1-\frac{j}{n} \biggr)\sum^j_{i=1}X_i- \frac{j}{n} \sum^n_{i=j+1}X_i.
\end{equation}
Here $c>0$ ($c=0.1$ in our simulation studies) is a small
number to avoid the boundary issue. For i.i.d. data, $j(1-j/n)$ is
proportional to the variance of $S_X(j)$, so $T^1_n$ is a
studentized version of $T^2_n$. For i.i.d. Gaussian data, $T^1_n$
is equivalent to likelihood ratio test; see Cs\"org\H{o} and Horv\'ath~\cite{csorgo}. Assume that,
under null hypothesis,
%
\begin{equation}\label{eq:cusuma}
\Biggl\{n^{-1/2}\sum^{\lfloor nt \rfloor}_{i=1} [ X_i - E(X_i) ] \Biggr\}
_{0\le t\le1}
\Rightarrow\tau\{B_t\}_{0\le t\le1},\qquad  \mbox{in the Skorohod space}
\end{equation}
for a standard Brownian motion $\{B_t\}_{t\ge0}$. The above convergence
requires finite-dimensional convergence and tightness; see Billingsley~\cite{bill}.
By the continuous mapping theorem,
$T^1_n\Rightarrow\max_{c\le t\le1-c}|B_t-tB_1|/\sqrt{t(1-t)}$
and $T^2_n/\sqrt{n}\Rightarrow\max_{c\le t\le1-c}|B_t-tB_1|$.

For the modulated stationary case (\ref{eq:null}), (\ref{eq:cusuma}) is
no longer valid. Moreover, since $T^1_n$ and $T^2_n$ do not
take into account the time-dependent variances $\sigma^2_i$, an
abrupt change in variances may lead to a false rejection of
$H_0$ when the mean remains constant. For example, our
simulation study in Section~\ref{sec:power} shows that the
empirical false rejection probability for $T^1_n$ and $T^2_n$
is about $10\%$ for nominal level $5\%$. To alleviate the issue
of non-constant variances, we adopt the self-normalization
approach as in previous sections. Recall $F_j$ and $\underline{V}_j$ in
(\ref{eq:Fj}). For each fixed $cn\le j\le(1-c)n$, by Theorem
\ref{thm:0} and Slutsky's theorem,
$F_j/\underline{V}_j\Rightarrow N(0,\tau^2)$ in distribution,
assuming the negligibility of the approximation errors.
Therefore, the self-normalization term $\underline{V}_j$ can
remove the time-dependent variances. In light of this, we can
simultaneously self-normalize the two terms $\sum^j_{i=1}X_i$
and $\sum^n_{i=j+1}X_i$ in (\ref{eq:SXj}) and propose the
self-normalized test statistic
%
\begin{equation}\label{eq:Tstar}
T^\mathrm{ SN}_n = \max_{c n \le j\le(1-c) n} \hat\tau^{-1}|T_n(j)|,
\quad\mbox{where }   T_n(j) = \frac{S_X(j)}{\sqrt{(1-j/n)^2 \underline{V}^2_j
+ (j/n)^2 \overline{V}^2_j}}.
\end{equation}
Here, $\underline{V}^2_j$ is defined as in (\ref{eq:Fj}),
$\overline{V}^2_j=\sum^n_{i=j+1}(X_i-\overline{X}_j)^2$ with
$\overline{X}_j=(n-j)^{-1}\sum^n_{i=j+1}X_i$.

\begin{theorem}\label{thm:test}
Assume (\ref{eq:sip}) holds. Let $\delta_n\to0$ be as in
(\ref{eq:thmtestcon}). Under $H_0$, we have
\[
\max_{cn\le j\le(1-c)n} |T_n(j)- \tau\widetilde{T}_n(j)|=\mathrm{O}_\mathrm{p}(\delta_n),
\]
where
\[
\widetilde{T}_n(j) = \frac{(1-j/n) \sum
^j_{i=1}\sigma_i (B_i-B_{i-1})
- j/n \sum^n_{i=j+1} \sigma_i (B_i-B_{i-1})}
{\sqrt{(1-j/n)^2\sum^j_{i=1} \sigma^2_i + (j/n)^2\sum^n_{i=j+1} \sigma^2_i}}.
\]
\end{theorem}

By Theorem~\ref{thm:test}, under $H_0$, $T^\mathrm{ SN}_n$ is
asymptotically equivalent to $\max_{cn\le j\le
(1-c)n}|\widetilde{T}_n(j)|$.\vspace*{1pt} Due to the self-normalization, for each
$j$, the time-dependent variances are removed and $\widetilde
{T}_n(j)\sim
N(0,1)$ has a standard normal distribution. However, $\widetilde{T}_n(j)$
and $\widetilde{T}_n(j')$ are correlated for $j\ne j'$. Therefore,
$\{\widetilde{T}_n(j)\}$ is a non-stationary Gaussian process with a
standard normal marginal density. Due to the large number of
unknown parameters $\sigma_i$, it is infeasible to obtain the
null distribution directly. On the other hand, Theorem
\ref{thm:test} establishes the fact that, asymptotically, the
distribution of $T^\mathrm{ SN}_n$ in (\ref{eq:Tstar}) depends only
on $\sigma_1,\ldots,\sigma_n$ and is robust against the
dependence structure of $\{e_i\}$, which motivates us to use
the wild bootstrap method in Section~\ref{sec:wild} to find the
critical value of $T^\mathrm{ SN}_n$.
\begin{enumerate}[(iii)]
\item[(i)] Compute $T_n(j)$ and find $\hat{J}=\argmax_{cn\le
j\le(1-c)n} |T_n(j)|$.
\item[(ii)] Divide the data into two blocks
$X_1,\ldots,X_{\hat{J}}$ and $X_{\hat{J}+1},\ldots,X_n$. Within
each block, subtract the sample mean from the
observations therein to obtain centered data. Pool all
centered data together and denote them by
$\epsilon_1,\ldots,\epsilon_n$.
\item[(iii)] Based on $\epsilon_1,\ldots,\epsilon_n$,
obtain an estimate $\hat\tau$ of $\tau$. See Section
\ref{sec:lrv} below.
\item[(iv)] Compute the test statistic $T^\mathrm{ SN}_n$ in
(\ref{eq:Tstar}).
\item[(v)] Based on $\epsilon_i$ in (ii), use the wild
bootstrap method in Section~\ref{sec:wild} to generate
synthetic data $\xi_1,\ldots,\xi_n$, and use (i)--(iv)
to compute the bootstrap test statistic $T^b_n$ based
on the bootstrap data $\xi_1,
\ldots,\xi_n$.
\item[(vi)] Repeat (v) many times and find
$(1-\alpha)$ quantile of those $T^b_n$s.
\end{enumerate}

As argued in Section~\ref{sec:wild}, the synthetic
data-generating scheme (v) inherits the time-varying
non-stationarity structure of the original data. Also, the
statistic $T^\mathrm{ SN}_n$ is robust against the dependence
structure, which justifies the proposed bootstrap method. If
$H_0$ is rejected, the change point is then estimated by
$\hat{J}=\argmax_{cn \le
j\le(1-c)n} |T_n(j)|$.

If there is no evidence to reject $H_0$, we briefly discuss how
to apply the same methodology to test $\tilde{H}_0\dvt  \sigma_1=\cdots
= \sigma_J \ne
\sigma_{J+1}=\cdots=\sigma_n$, that is, whether there is a
change point in the variances $\sigma^2_i$. By
(\ref{eq:xinons}), we have $(X_i-\mu)^2=\sigma^2_i + \sigma^2_i
\zeta_i$, where $\zeta_i=e^2_i-1$ has mean zero. Therefore,
testing a change point in the variances $\sigma^2_i$ of $X_i$
is equivalent to testing a change point in the mean of the new
data $\tilde{X}_i=(X_i-\bar{X})^2$.

\subsection{Long-run variance estimation}\label{sec:lrv}

To apply the results in Sections
\ref{sec:cltx}--\ref{sec:cusum}, we need a consistent estimate
of the long-run variance~$\tau^2$. Most existing works deal with
stationary time
series through various block bootstrap and subsampling
approaches; see Lahiri~\cite{lahiri2} and references therein.
Assuming a near-epoch dependent mixing condition, De Jong and Davidson~\cite{dejong}
established the consistency of a kernel estimator of $\operatorname{ Var}(\sum
^n_{i=1}X_i)$, and their result can be used to
estimate $\tau^2_n$ in (\ref{eq:lrvnon}) for the CLT of $\sqrt{n}(\bar
{X}-\mu)$. However,
for the change point problem in Section~\ref{sec:cusum}, we need an
estimator of
the long-run variance $\tau^2$ of the unobservable process $\{e_i\}$,
so the method in De Jong and Davidson~\cite{dejong} is not
directly applicable.

To
attenuate the non-stationarity issue, we extend the idea in
Section~\ref{sec:cltx} to blockwise self-normalization. Let
$k_n$ be the block length. Denote by $\ell_n=\lfloor
n/k_n\rfloor$ the largest integer not exceeding $n/k_n$. Ignore
the boundary and divide $1,\ldots,n$ into $\ell_n$ blocks
%
\begin{equation}\label{eq:Ij}
\mathcal{I}_j=\{(j-1)k_n+1,\ldots,jk_n\},\qquad   j=1,\ldots,\ell_n.
\end{equation}
Recall the overall sample mean $\bar{X}$. For each block $j$,
define the self-normalized statistic
%
\begin{equation}\label{eq:Dj}
D_j=\frac{k_n [ \bar{X}(j) - \bar{X} ]}{V(j)},
 \qquad\mbox{where }
\bar{X}(j)=\frac{1}{k_n} \sum_{i\in\mathcal{I}_j} X_i,
V^2(j) = \sum_{i\in\mathcal{I}_j} [X_i-\bar{X}(j)]^2.
\end{equation}
By Proposition~\ref{cor:1}, the self-normalized statistics
$D_1, \ldots, D_{\ell_n}\sim N(0,\tau^2)$ are asymptotically i.i.d.
Thus, we propose
estimating $\tau^2$ by
%
\begin{equation}\label{eq:lrvesta}
\hat\tau^2= \frac{1}{\ell_n} \sum^{\ell_n}_{j=1} D^2_j.
\end{equation}
As in (\ref{eq:volwei})--(\ref{eq:Omegan}), we define the
quantities on block $j$
%
\begin{eqnarray}\label{eq:volweij}
r(j) &=& |\sigma_{jk_n}| +
\sum_{i\in\mathcal{I}_j} |\sigma_i-\sigma_{i-1}|  \quad \mbox{and}
\quad
r^*(j) = |\sigma^2_{jk_n}| +
\sum_{i\in\mathcal{I}_j} |\sigma^2_i-\sigma^2_{i-1}|,
\\
\label{eq:Sigmaj}
\Sigma^2(j)&=& \sum_{i\in\mathcal{I}_j} \sigma^2_i   \quad\mbox{and}\quad
\Sigma^{*2}(j)=\biggl( \sum_{i\in\mathcal{I}_j}\sigma^4_i \biggr)^{1/2}.
\end{eqnarray}

\begin{theorem}\label{thm:3}
Let (\ref{eq:sip}) hold with $\Delta_n= n^{1/4}\log(n)$. Recall
$r_n, \Sigma_n$ in (\ref{eq:volwei})--(\ref{eq:Omegan}). Define
%
\begin{equation}\label{eq:Mn}
M_n = \frac{1}{k_n} + \max_{1\le j\le\ell_n} \frac{\Sigma^{*2}(j)
+r^*(j) \Delta_n}{\Sigma^2(j)}
+ \max_{1\le j\le\ell_n} \frac{r(j)\Delta_n}{\Sigma(j)}.
\end{equation}
Assume that $r_n\Delta_n/\Sigma_n\to0$ and
%
\begin{equation}
\chi_n=\ell_n^{-1/2}
+ \log(n) M_n + \sqrt{\log(n)} \frac{\Sigma_n}{\ell^{2}_n}
\sum^{\ell_n}_{j=1}\frac{1}{\Sigma(j)} + \frac{\Sigma^2_n}{\ell^{3}_n}
\sum^{\ell_n}_{j=1}\frac{1}{\Sigma^2(j)} \to0.
\end{equation}
Then $\hat\tau^2-\tau^2 = \mathrm{O}_\mathrm{ p}(\chi_n)$. Consequently,
$\hat\tau$ is a consistent estimate of $\tau$.
\end{theorem}

Consider Example~\ref{exmp:6} with $\gamma\in[0,1)$. Then
$\chi_n\asymp
\sqrt{\log(n)/\ell_n}+\log^{2}(n)(n^{1/4}/\sqrt{k_n}+n^{5/4-\gamma
}/k_n+\sqrt{k_n}
n^{1/4-\gamma})$. For $\gamma\in(3/4,1)$, it can be shown that
the optimal rate is $\chi_n\asymp n^{-1/8}\log^{5/4}(n)$ when
$k_n\asymp n^{3/4}\log^{3/2}(n)$. In Example
\ref{exmp:4} with $\sigma_i=i^\beta$ for some $\beta\in[0,1)$,
elementary but tedious calculations show that the optimal
rate is
\[
\chi_n \asymp\cases{
n^{-1/8} \log^{5/4}(n), \qquad
k_n\asymp n^{3/4}\log^{3/2}(n), \vspace*{2pt}\cr
\quad \beta\in[0,3/4],\vspace*{5pt}\cr
n^{{(\beta-1)}/{(5-4\beta)}}\{\log(n)\}^{{(8(1-\beta))}/{(5-4\beta)}}, \qquad
k_n\asymp n^{{(4.5-4\beta)}/{(5-4\beta)}}\{\log(n)\}^{{4}/{(5-4\beta
)}}, \vspace*{2pt}\cr
\quad \beta\in(3/4,1).}
\]

\subsection{Some possible extensions}\label{sec:ext}\vspace*{-2pt}
The self-normalization approaches in Sections
\ref{sec:cltx}--\ref{sec:lrv} can be extended to linear
regression model (\ref{eq:lr}) with modulated stationary time series
errors. The approach in Phillips, Sun and Jin~\cite{phi} is not
applicable here due to non-stationarity. For simplicity, we
consider the simple case that $p=2, U_i=(1,i/n), $ and
$\beta=(\beta_0,\beta_1)'$. Hansen~\cite{hansen1995} studied a similar setting
for martingale difference errors.
Denote by $\hat\beta_0$ and
$\hat\beta_1$ the simple linear regression estimates of
$\beta_0$ and $\beta_1$ given by
%
\begin{equation}\label{eq:lrest}
\hat\beta_1 = \frac{n\sum^n_{i=1} i X_i - \sum^n_{i=1} i \sum^n_{i=1}X_i}
{\sum^n_{i=1} i^2 - (\sum^n_{i=1} i)^2/n} \quad \mbox{and}\quad
\hat\beta_0 = \bar{X}-\hat\beta_1 (n+1)/(2n).
\end{equation}
Then simple algebra shows that
\[
\hat\beta_0 - \beta_0 = \frac{2}{n^2-n} \sum^n_{i=1} (2n-3i+1) \sigma_i
e_i,\qquad
\hat\beta_1 - \beta_1 = \frac{6}{n^2-1} \sum^n_{i=1} (2i-n-1) \sigma_i e_i.
\]
The latter expressions are linear combinations of
$\{e_i\}$. Thus, by the same argument in Proposition~\ref{cor:1}
and Theorem~\ref{thm:0}, we have self-normalized CLTs for
$\hat\beta_0$ and $\hat\beta_1$.

\begin{theorem}\label{thm:lr}
Let $s_{i,0}=(2n-3i+1) \sigma_i$ and $s_{i,1}=(2i-n-1)
\sigma_i$. Assume that $\{s_{i,0}\}_{1\le i\le n}$ and
$\{s_{i,1}\}_{1\le i\le n}$ satisfy condition
(\ref{eq:thmtestcon}). Then as $n\to\infty$,
\begin{eqnarray*}
 \frac{n^2(\hat\beta_0-\beta_0)}{2 V_{n,0}} &\Rightarrow & N(0,\tau^2),
\qquad \mbox{where }   V_{n,0}^2 = \sum^n_{i=1}
(2n-3i+1)^2 (X_i-\hat\beta_0-\hat\beta_1 i/n)^2, \\[-3pt]
 \frac{n^2(\hat\beta_1-\beta_1)}{6 V_{n,1}}& \Rightarrow& N(0,\tau^2),
 \qquad\mbox{where }   V_{n,1}^2 = \sum^n_{i=1}
(2i-n-1)^2 (X_i-\hat\beta_0-\hat\beta_1 i/n)^2.
\end{eqnarray*}
\end{theorem}

The long-run variance $\tau^2$ can be estimated using the idea of
blockwise self-nor\-malization in Section~\ref{sec:lrv}. Let
$k_n, \ell_n$ and $\mathcal{I}_j$ be defined as in Section
\ref{sec:lrv}. Then we propose
%
\begin{equation}\label{eq:lrvest}
\hat\tau^2= \frac{1}{\ell_n} \sum^{\ell_n}_{j=1} D^2_j,
\qquad \mbox{where }   D_j = \frac{\sum_{i\in\mathcal{I}_j} (X_i-\hat
\beta_0-\hat\beta_1 i/n)}
{\sqrt{\sum_{i\in\mathcal{I}_j} (X_i-\hat\beta_0-\hat\beta_1 i/n)^2}}.
\end{equation}
Here, $D_1,\ldots,D_{\ell_n}$ are asymptotically i.i.d. normal
random variables with mean zero and variance~$\tau^2$.
Consistency can be established under similar conditions as in
Theorem~\ref{thm:3}.

For the general linear regression model (\ref{eq:lr}),
the linearly weighted average structure of linear regression
estimates allows us to obtain self-normalized CLTs as
in Theorem~\ref{thm:lr} under more complicated conditions.
Also, it is possible to extend the proposed method to the
nonparametric regression model with time-varying variances
%
\begin{equation}\label{eq:npm}
X_i=f(i/n)+\sigma_i e_i,
\end{equation}
where $f(\cdot)$ is a nonparametric time trend of interest.
Nonparametric estimates, for example, the Nadaraya--Watson\vadjust{\goodbreak}
estimate, are usually based on locally weighted observations.
The latter feature allows us to derive similar self-normalized~CLT. However, the change point problem for~(\ref{eq:lr}) and
(\ref{eq:npm}) will be more challenging, and
Aue \textit{et al.}~\cite{aue2008a} studied~(\ref{eq:lr}) for uncorrelated errors with
constant variance.
Also, it is more difficult to address the bandwidth selection issues;
see Altman~\cite{altman} for related
contribution when $\sigma_i\equiv1$.
It remains a direction
of future research to investigate (\ref{eq:lr}) and
(\ref{eq:npm}).\vspace*{-1pt}

\section{Simulation study}\label{sec:simu}\vspace*{-1pt}

\subsection{\texorpdfstring{Selection of block length $k_n$ for $\hat\tau$}{Selection of block length k n for tau}}\label{sec:kn}\vspace*{-1pt}
Recall that
$D_1,\ldots,D_{\ell_n}$ in (\ref{eq:lrvest}) are asymptotically
i.i.d. normal random variables. To get a sensible choice of the block
length parameter $k_n$,
we propose a simulation-based method by minimizing the
empirical mean squared error (MSE):
\begin{enumerate}[(iii)]
\item[(i)] Simulate $n$ i.i.d. standard normal random
variables $Z_1,\ldots,Z_n$.
\item[(ii)] Based on $Z_1,\ldots,Z_n$, obtain $\hat\tau$
with block length $k$.
\item[(iii)] Repeat (i)--(ii) many times,
compute empirical $\operatorname{MSE}(k)$ as the average of
realizations of $(\hat\tau-1)^2$, and find the optimal
$k$ by minimizing $\operatorname{MSE}(k)$.
\end{enumerate}
We find that the optimal block length $k$ is about 12 for
$n=120$, about 15 for $n=240$, about 20 for $n=360,600$ and
about 25 for $n=1200$.\vspace*{-1pt}

\subsection{Empirical coverage
probabilities}\label{sec:ecp}\vspace*{-1pt}
Let sample size $n=120$. Recall $e_i$ and $\sigma_i$ in
(\ref{eq:xinons}). For $\sigma_i$, consider four choices:
\begin{eqnarray*}
&&\mathrm{A1}\!:\quad  \sigma_i=0.2 \mathbf{1}_{i\le n/2} + 0.6
\mathbf{1}_{i>n/2},\\
&&\mathrm{A2}\!:\quad  \sigma_i=0.2\{1+\cos^2(i/n^{4/5})\}, \\
&&\mathrm{A3}\!:\quad \sigma_i=0.2+0.1\log(1+|i-n/2|), \\
&&\mathrm{A4}\!:\quad \sigma_i=0.3+\phi(i/60),
\end{eqnarray*}
where $\phi$ is the standard normal density, and $\mathbf{1}$ is
the indicator function. The sequences A1--A4 exhibit different
patterns, with a piecewise constancy for A1, a cosine shape for
A2, a sharp
change around time $n/2$ for A3 and a gradual downtrend for A4.
Let $\varepsilon_i$ be i.i.d.
N(0, 1). For $e_i$, we consider both linear and nonlinear processes.
\begin{eqnarray*}
&&\mathrm{ B1}\!:\quad  e_i=\{\eta_i-E(\eta_i)\}/\sqrt{\operatorname{ Var}(\eta_i)},
\qquad\mbox{where }
\eta_i = \theta|\eta_{i-1}| + \sqrt{1-\theta^2} \varepsilon_i, |\theta
|<1.\\
&&\mathrm{B2}\!:\quad e_i = \sum^\infty_{j=0} a_j \varepsilon_{i-j},
 \qquad\mbox{where }   a_j= \frac{(j+1)^{-\beta}}{\sqrt{\sum^\infty
_{j=0} (j+1)^{-2\beta}}}, \beta>1/2.
\end{eqnarray*}
For B1, by Wu~\cite{wuw}, (\ref{eq:pro1con}) holds. By Andel, Netuka and Svara~\cite{andel},
$E(\eta_i)=\theta\sqrt{2/\uppi}$ and $\operatorname{
Var}(\eta_i)=1-2\theta^2/\uppi$. To examine how the strength of
dependence affects the performance, we consider $\theta=0, 0.4,
0.8$, representing independence, intermediate and strong
dependence, respectively.\vadjust{\goodbreak} For B2 with $\beta>2$, (\ref{eq:sip})
holds with $\Delta_n=n^{1/4}\log(n)$, and we consider three
cases $\beta=2.1, 3, 4$. To assess the effect of block length
$k_n$, three choices $k_n=8,10,12$ are used. Thus, we consider
all 72 combinations of $\{\mathrm{ A1},\mathrm{ A2},\mathrm{ A3}, \mathrm{
A4}\} \times\{\mathrm{ B1}, \theta=0,0.4,0.8; \mathrm{ B2},
\beta=2.1,3,4\} \times\{k_n=8,10,12\}$.

\begin{table}
\caption{Coverage probabilities (in percentage) for $\mu$ in
(\protect\ref{eq:xinons}) with $e_i$ from B1 [(a)] and
B2 [(b)]. Nominal level is $95\%$.
SN and WB denote self-normalization-based confidence intervals
using asymptotic theory in Proposition \protect\ref{cor:1} and the
wild bootstrap procedure, respectively; ST, BB, SBB denote
stationarity-based confidence intervals using asymptotic theory
in (\protect\ref{eq:fanyao}), non-overlapping block bootstrap and
studentized non-overlapping block bootstrap, respectively}\label{tab1}
\fontsize{8}{10}{\selectfont{
\begin{tabular*}{\textwidth}{@{\extracolsep{\fill}}llllllllllllll@{}}
\hline
$\theta$ &$k_n$ &$\sigma_i$ &SN &WB &ST &BB &SBB &$\sigma_i$ &SN &WB
&ST &BB &SBB \\
\hline
\multicolumn{14}{c@{}}{(a) Model B1}\\
0.0 &\phantom{0}8 & A1&98.0 &94.7 &93.1 &92.2 &92.8 &A2 &96.6 &95.2 &92.3 &92.5 &92.5
\\
&10 & &98.2 &95.0 &92.6 &92.4 &92.2 & &94.6 &94.6 &90.0 &89.5
&89.4 \\
&12 & &98.1 &95.6 &91.7 &91.4 &91.1 & &92.1 &93.7 &89.7 &89.5 &89.6
\\[2pt]
&\phantom{0}8 &A3 &96.4 &95.0 &92.5 &92.3 &92.0 &A4 &96.6 &95.6 &93.1 &92.6 &93.0\\
&10 & &94.7 &94.7 &90.8 &90.6 &90.6 & &95.1 &95.1 &91.4 &91.3
&91.3\\
&12 & &93.7 &94.8 &90.8 &90.4 &90.5 & &92.9 &93.7 &89.8 &89.7 &89.5\\
[4.5pt]
0.4 &\phantom{0}8 & A1&98.7 &95.9 &92.7 &92.6 &92.9 & A2&96.6 &95.3 &92.5 &92.4
&92.0\\
&10 & &98.5 &95.7 &92.8 &92.7 &92.3 & &95.4 &95.4 &91.6 &91.1
&91.6\\
&12 & &98.0 &95.0 &90.8 &90.8 &90.2 & &92.5 &94.0 &89.4 &89.1 &89.4\\
[2pt]
&\phantom{0}8 &A3 &96.6 &95.2 &91.7 &91.7 &91.6 &A4 &95.4 &94.1 &90.8 &90.9 &90.6\\
&10 & &95.3 &95.5 &91.5 &91.3 &91.5 & &95.0 &94.8 &91.2 &90.7
&90.8\\
&12 & &93.1 &94.6 &90.2 &89.9 &89.9 & &94.1 &95.1 &90.3 &89.8 &90.1\\
[4.5pt]
0.8 &\phantom{0}8 & A1&97.9 &94.6 &87.8 &86.8 &87.3 &A2 &96.1 &94.7 &87.2 &87.3
&87.0\\
&10 & &97.6 &95.5 &87.3 &87.0 &86.7 & &93.3 &92.9 &86.4 &86.8
&86.1\\
&12 & &97.3 &94.0 &85.8 &85.5 &85.1 & &92.6 &93.4 &86.5 &86.4 &86.4 \\
[2pt]
&\phantom{0}8 &A3 &94.8 &93.5 &85.7 &85.7 &86.0 & A4&95.5 &94.7 &86.3 &86.1 &86.1 \\
&10 & &93.5 &93.8 &85.7 &85.5 &85.2 & &95.3 &95.1 &88.5 &88.3
&88.5 \\
&12 & &92.4 &93.3 &87.2 &86.7 &86.9 & &92.6 &94.2 &86.3 &85.8 &85.7 \\
\end{tabular*}}}\vspace*{7pt}
\fontsize{8}{10}{\selectfont{
\begin{tabular*}{\textwidth}{@{\extracolsep{\fill}}llllllllllllll@{}}
\multicolumn{14}{@{}l}{$\beta$}\\
\hline
\multicolumn{14}{c}{(b) Model B2}\\
4.0 &\phantom{0}8 & A1&97.6 &94.9 &91.8 &91.4 &91.9 &A2 &95.9 &94.2 &91.9 &92.0 &91.1
\\
&10 & &97.7 &93.2 &88.9 &88.1 &88.3 &&95.7 &95.7 &92.1 &91.8 &92.1
\\
&12 & &97.9 &95.5 &90.7 &90.2 &90.0 & &93.3 &94.6 &90.0 &89.9 &89.7 \\
[2pt]
&\phantom{0}8 &A3 &94.6 &93.3 &89.8 &89.5 &89.5 & A4&95.6 &94.7 &91.3 &91.7 &91.0 \\
&10 & &95.1 &95.2 &91.6 &91.4 &91.5 &&95.4 &95.9 &92.8 &92.2 &93.0
\\
&12 & &93.8 &95.4 &90.8 &90.6 &90.2 & &93.9 &94.9 &88.9 &88.5 &88.6 \\
[4.5pt]
3.0 &\phantom{0}8 &A1 &99.1 &95.7 &91.1 &91.0 &91.2 & A2&95.8 &94.6 &90.4 &89.8
&90.1\\
&10 & &98.5 &96.4 &91.6 &90.9 &91.1 & &95.6 &95.2 &92.1 &91.9
&91.5\\
&12 & &97.9 &94.6 &89.6 &89.3 &89.0 & &94.1 &95.0 &90.5 &90.2 &90.4\\
[2pt]
&\phantom{0}8 & A3&95.9 &94.6 &92.0 &91.9 &91.7 & A4&96.0 &94.5 &90.6 &90.4 &90.3\\
&10 & &94.3 &94.4 &90.0 &89.9 &89.8 & &94.3 &94.4 &89.2 &89.3
&88.9\\
&12 & &93.2 &94.5 &88.9 &88.6 &88.7 & &93.1 &94.1 &89.6 &88.9 &88.8\\
[4.5pt]
2.1 &\phantom{0}8 & A1&97.1 &92.5 &86.2 &86.2 &85.5 & A2&95.7 &93.8 &88.9 &89.0 &88.7
\\
&10 & &97.6 &94.7 &89.2 &88.9 &88.6 &&93.5 &93.6 &88.8 &88.8 &88.4
\\
&12 & &97.2 &95.1 &87.9 &87.5 &87.7 & &92.6 &93.9 &88.0 &87.6 &87.7 \\
[2pt]
&\phantom{0}8 & A3&94.0 &93.7 &88.5 &88.4 &88.3 &A4 &95.0 &93.1 &88.8 &88.7 &88.6 \\
&10 & &93.3 &93.8 &88.1 &87.9 &87.8 & &94.1 &94.2 &89.1 &88.8
&89.1\\
&12 & &92.9 &94.7 &89.1 &88.4 &88.4 & &91.5 &92.6 &87.7 &87.5 &87.5 \\
\hline
\end{tabular*}}}
\end{table}

Without loss of generality we examine coverage probabilities
based on $10^3$ realized confidence intervals for $\mu=0$ in
(\ref{eq:xinons}). We compare our self-normalization-based
confidence intervals to some stationarity-based methods. For
(\ref{eq:xinons}), if we pretend that the error process
$\{\tilde{e}_i=\sigma_i e_i\}$ is stationary, then we can use
(\ref{eq:fanyao}) to construct an asymptotic confidence interval
for $\mu$. Under stationarity, the long-run variance $\tau^2$
of $\{\tilde{e}_i\}$ can be similarly estimated through the
block method in Section~\ref{sec:lrv} by using the
non-normalized version $D_j=\sqrt{k_n} [ \bar{X}(j) - \bar{X}
]$ in (\ref{eq:lrvest}); see Lahiri~\cite{lahiri2}. Thus, we compare
two self-normalization-based methods and three stationarity-based alternatives:
self-normalization-based confidence
intervals through the asymptotic theory in Proposition
\ref{cor:1} (SN) and the wild bootstrap (WB) in Section
\ref{sec:wild}; stationarity-based confidence intervals through
the asymptotic theory (\ref{eq:fanyao}) (ST), non-overlapping
block bootstrap (BB) and studentized non-overlapping block
bootstrap (SBB) in Section~\ref{sec:wild}. From the results in
Table~\ref{tab1}, we see that the coverage
probabilities of the proposed self-normalization-based methods
(columns SN and WB) are close to the nominal level $95\%$ for
almost all cases considered. By contrast, the stationarity-based
methods (columns ST, BB and SBB) suffer from substantial
undercoverage, especially when dependence is strong
($\theta=0.8$ in Table~\ref{tab1}(a) and $\beta=2.1$ in Table~\ref{tab1}(b)). For the two self-normalization-based methods, WB
slightly outperforms SN.

\begin{table}
\caption{Size (in percentage) comparison of $T^1_n$ and $T^2_n$
in (\protect\ref{eq:cusum}) and $T^\mathrm{ SN}_n$ in (\protect\ref{eq:Tstar}),
with sample size $n=120$, nominal level $5\%$, and block length
$k_n=10$}\label{tab:3}
\begin{tabular*}{\textwidth}{@{\extracolsep{\fill}}llld{2.1}d{2.1}lld{2.1}d{2.1}@{}}
\hline
&\multicolumn{4}{l}{Model B1}
&\multicolumn{4}{l@{}}{Model B2} \\[-6pt]
&\multicolumn{4}{c}{\hrulefill}
&\multicolumn{4}{c@{}}{\hrulefill} \\
\multicolumn{1}{@{}l}{$\sigma_i$} & \multicolumn{1}{l}{$\theta$} &
\multicolumn{1}{l}{$T^{\mathrm{SN}}_n$} &
\multicolumn{1}{l}{$T^{1}_n$} &
\multicolumn{1}{l}{$T^2_n$} &
\multicolumn{1}{l}{$\beta$} &
\multicolumn{1}{l}{$T^{\mathrm{SN}}_n$} &
\multicolumn{1}{l}{$T^{1}_n$} & \multicolumn{1}{l@{}}{$T^2_n$} \\
\hline
A1&0.0 &4.9 &9.1 &8.4 &2.1 &7.3 &12.2 &13.4 \\
 &0.4 &4.7 &9.4 &9.6 &3.0 &4.7 &8.6 &9.2 \\
&0.8 &6.0 &15.1 &14.7 &4.0 &5.6 &9.9 &7.7 \\[3pt]
A2&0.0 &5.7 &8.2 &6.1 &2.1 &5.8 &9.5 &8.6 \\
 &0.4 &6.1 &8.9 &6.8 &3.0 &5.3 &9.6 &6.8 \\
&0.8 &7.3 &12.6 &9.3 &4.0 &4.2 &7.5 &4.2 \\[3pt]
A3&0.0 &5.0 &5.7 &4.8 &2.1 &5.5 &7.7 &6.7 \\
 &0.4 &5.3 &6.9 &5.4 &3.0 &5.8 &6.1 &4.9 \\
&0.8 &7.0 &9.8 &10.0 &4.0 &5.0 &6.5 &4.2 \\[3pt]
A4&0.0 &5.4 &8.4 &6.0 &2.1 &6.9 &8.8 &7.1 \\
 &0.4 &5.7 &7.9 &5.2 &3.0 &4.8 &6.6 &6.3 \\
&0.8 &7.2 &11.1 &9.2 &4.0 &5.3 &6.2 &5.8 \\
\hline
\end{tabular*}
\end{table}

\subsection{Size and power study}\label{sec:power}

\begin{figure}

\includegraphics{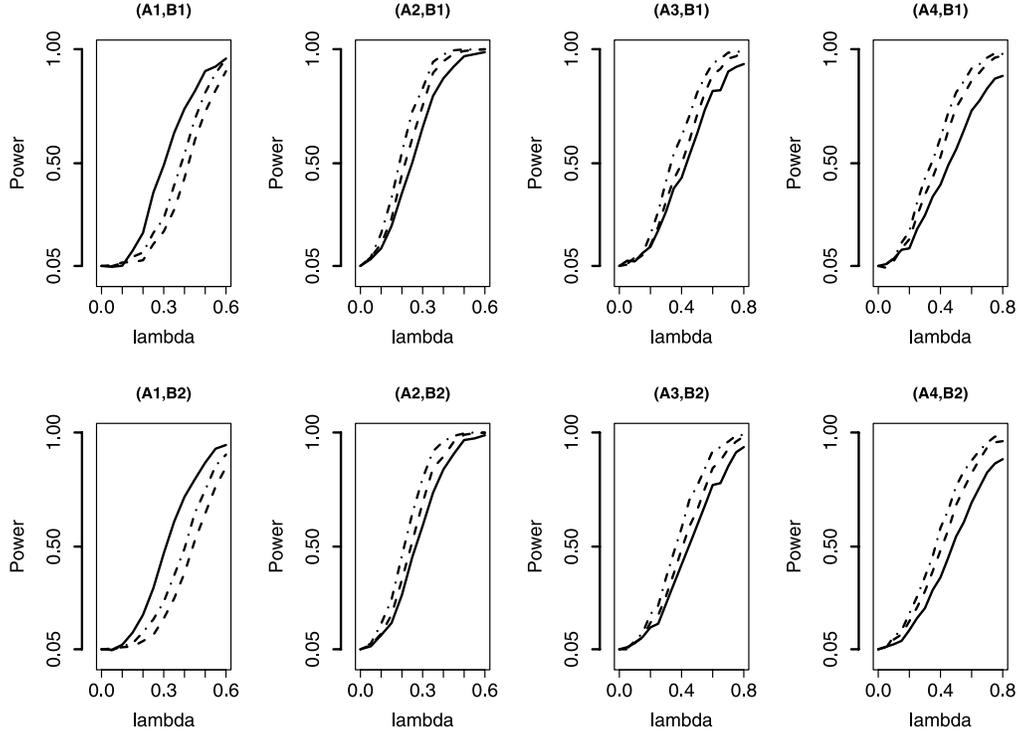}

\caption{Size-adjusted power curves for $T^1_n$ (dashed curve) and
$T^2_n$ (dotdash curve) in (\protect\ref{eq:cusum})
and $T^\mathrm{ SN}_n$ (solid curve) in (\protect\ref{eq:Tstar}) as functions of
change size $\lambda$ (horizontal axis) with
sample size $n=120$ and block length
$k_n=10$. For (A1, B1)--(A4, B1), the error process $\{e_i\}$ is from B1
with $\theta=0.4$;
for (A1, B2)--(A4, B2), the error process $\{e_i\}$ is from B2 with
$\beta=3.0$.}\label{fig:power}
\end{figure}

In (\ref{eq:null}), we use the same setting for $\sigma_i$
and $e_i$ as in Section~\ref{sec:ecp}. For mean $\mu_i$, we
consider $\mu_i = \lambda\mathbf{1}_{i>40},\lambda\ge0$, and
compare the test statistics $T^1_n, T^2_n$ in (\ref{eq:cusum})
and $T^\mathrm{ SN}_n$ in (\ref{eq:Tstar}).\vspace*{1pt} First, we compare their
size under the null with $\lambda=0$. The
critical value of $T^\mathrm{ SN}_n$ is obtained using the wild bootstrap
in Section
\ref{sec:cusum}; for $T^1_n$ and $T^2_n$, their critical values
are based on the block bootstrap in Section~\ref{sec:wild}. In
each case, we use $10^3$ bootstrap samples, nominal level
$5\%$, and block length $k_n=10$, and summarize the empirical
sizes (under the null $\lambda=0$) in Table~\ref{tab:3} based
on $10^3$ realizations. While $T^\mathrm{ SN}_n$ has size close to
$5\%$, $T^1_n$ and $T^2_n$ tend to over-reject the null, and the
false rejection probabilities can be three times the
nominal level of $5\%$. Next, we compare the size-adjusted power.
Instead of using the bootstrap methods to obtain critical
values, we use $95\%$ quantiles of $10^4$ realizations of the
test statistics when data are simulated directly from the null
model so that the empirical size is exactly $5\%$. Figure
\ref{fig:power} presents the power curves for combinations
\{A1--A4\} $\times$ \{B1 with $\theta=0.4$; B2 with
$\beta=3.0$\} with $10^3$ realizations each. For~A1, $T^\mathrm{
SN}_n$ outperforms $T^1_n$ and $T^2_n$; for A2--A4, there is a
moderate loss of power for $T^\mathrm{SN}_n$. Overall, $T^\mathrm{
SN}_n$ has power comparable to other two tests. In practice,
however, the null model is unknown, and when one turns to the
bootstrap method to obtain the critical values, the usual
CUSUM tests $T^1_n$ and $T^2_n$ will likely
over-reject the null as shown in Table~\ref{tab:3}. In
summary, with such small sample size and complicated
time-varying variances structure, $T^\mathrm{ SN}_n$ along with the wild
bootstrap method delivers reasonably good power and the size is
close to nominal level.

Finally, we point out that the proposed self-normalization-based
methods are not robust to models
with time-varying correlation structures. For example, consider the
model $e_i=0.3e_{i-1}+\varepsilon_i$ for $1\le i\le60$
and $e_i=0.8e_{i-1}+\varepsilon_i$ for $61\le i\le n$, where
$\varepsilon_i$ are i.i.d. N(0, 1).
With $k_n=10$, the sizes (nominal level $5\%$) for the three tests
$T^\mathrm{ SN}_n$, $T^{1}_n$, $T^2_n$ are
0.154, 0.196, 0.223 for A1. Future research directions include (i)
developing tests for change in the variance or covariance structure for
(\ref{eq:xinons}) (See Incl\'an and Tiao~\cite{inclan}, Aue \textit{et al.}~\cite{aue2009} and
Berkes, Gombay and Horv\'ath~\cite{berkes} for related contributions);
and (ii) developing methods that are robust to changes in
correlations.

\section{Applications to two real data sets}\label{sec:app}
\subsection{Annual mean precipitation in Seoul during 1771--2000}
\begin{figure}

\includegraphics{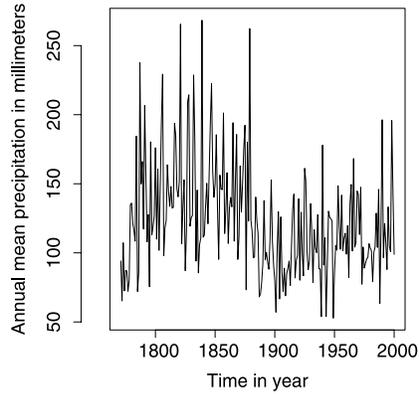}

\caption{Annual mean precipitation in Seoul from
1771--2000.}\label{fig:seoul}\vspace*{-3pt}
\end{figure}

The data set consists of annual mean precipitation rates in
Seoul during 1771--2000; see
Figure~\ref{fig:seoul} for a plot. The mean levels
seem to be different for the two time periods 1771--1880 and
1881--2000. Ha and Ha~\cite{ha} assumed the observations are
i.i.d. under the null
hypothesis. As shown in Figure~\ref{fig:seoul}, the variations change
over time. Also, the auto-correlation function plot (not reported here)
indicates strong dependence up to lag 18.
Therefore, it is more reasonable to apply our
self-normalization-based test that is tailored to deal with
modulated stationary processes. With sample size $n=230$, by the
method in Section~\ref{sec:kn}, the optimal block length is
about 15. Based on $10^5$ bootstrap samples as described in
Section~\ref{sec:cusum}, we obtain the corresponding p-values
0.016, 0.005, 0.045, 0.007, with block length
$k_n=12, 14, 16, 18$, respectively. For all choices of $k_n$,
there is compelling evidence that a change point occurred at
year 1880. While our result is consistent with that of Ha and Ha~\cite{ha},
our modulated stationary time series framework seems to be
more reasonable. Denote by $\mu_1$ and $\mu_2$ the mean levels
over pre-change and post-change time periods 1771--1880 and
1881--2000. For the two sub-periods with sample sizes 110 and
120, the optimal block length is about 12. With $k_n=12$,
applying the wild bootstrap in Section~\ref{sec:wild} with
$10^5$ bootstrap samples, we obtain $95\%$ confidence intervals
$[121.7,161.3]$ for $\mu_1$, $[100.9,114.3]$ for $\mu_2$. For
the difference $\mu_1-\mu_2$, with optimal block length
$k_n=15$, the $95\%$ wild bootstrap confidence interval is
$[19.6,48.2]$. Note that the latter confidence interval for
$\mu_1-\mu_2$ does not cover zero, which provides further
evidence for $\mu_1\ne\mu_2$ and the existence of a change point at
year 1880.

\begin{figure}[b]
\vspace*{-3pt}
\includegraphics{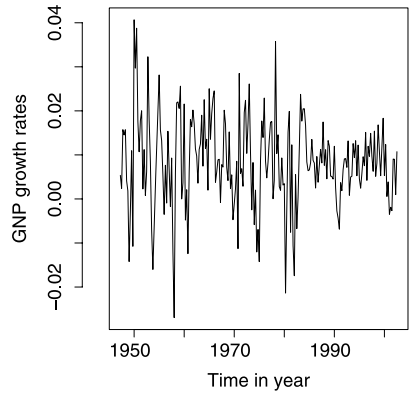}

\caption{Quarterly U.S. GNP growth rates from 1947--2002.}\label{fig:gnp}
\end{figure}

\subsection{Quarterly U.S. GNP growth rates during 1947--2002}\label{sec:gnp}

The data set consists of quarterly
U.S. Gross National Product (GNP) growth rates from the first
quarter of 1947 to the third quarter of 2002; see Section 3.8\vadjust{\goodbreak}
in Shumway and Stoffer~\cite{shumway} for a stationary autoregressive model approach.
However, the plot in
Figure~\ref{fig:gnp} suggests a non-stationary pattern: the
variation becomes smaller after year 1985 whereas the mean
level remains constant. Moreover, the stationarity test in
Kwiatkowski \textit{et al.}~\cite{kwi} provides fairly strong evidence for
non-stationarity with a p-value of 0.088. With the block length
$k_n=12,14,16,18$, we obtain the corresponding p-values
$0.853,0.922,0.903,0.782$, and hence there is no evidence to
reject the null hypothesis of a constant mean $\mu$. Based on
$k_n=15$, the $95\%$ wild bootstrap confidence interval for
$\mu$ is $[0.66\%,1.00\%]$. To test whether there is a
change point in the variance, by the discussion in the last
paragraph of Section~\ref{sec:cusum},\vadjust{\goodbreak} we consider
$\tilde{X}_i=(X_i-\underline{X}_n)^2$. With $k_n=12,14,16,18$,
the corresponding p-values are $0.001, 0.006,0.001,0.010$,
indicating strong evidence for a change point in the variance
at year 1984. In summary, we conclude that there is no change point
in the mean level, but there is a change point in the
variance at year 1984.

\begin{appendix}
\section*{Appendix: Proofs}\label{sec:proof}\vspace*{-5pt}

\begin{pf*}{Proof of Theorem \protect\ref{thm:0}}
Let $r_j = |\sigma_j| + \sum^j_{i=2} |\sigma_i-\sigma_{i-1}|$.
By the triangle inequality, we have $r_j \le r_n$. Recall $S_i$
in (\ref{eq:sip}). By the summation by parts
formula, (\ref{eq:thm0a}) follows via
%
\begin{eqnarray}\label{eq:tnsip}
F_j& =& \sum^j_{i=1} \sigma_i (S_i-S_{i-1})
=\sigma_j S_j +
\sum^{j-1}_{i=1} (\sigma_i - \sigma_{i+1}) S_i \nonumber\\
&=& \sigma_j \tau B_j + \sum^{j-1}_{i=1} (\sigma_i - \sigma_{i+1}) \tau
B_i +
\mathrm{O}_\mathrm{ a.s.}(r_n\Delta_n)\\
&= & \tau\sum^j_{i=1} \sigma_i (B_i-B_{i-1}) + \mathrm{O}_\mathrm{
a.s.}(r_n\Delta_n).\nonumber
\end{eqnarray}
By Kolmogorov's maximal inequality for independent random
variables, for $\delta>0$,
%
\begin{equation}\label{eq:maxin}
P \Biggl\{ \max_{1\le j\le n} \Biggl|\sum^j_{i=1} \sigma_i (B_i-B_{i-1})
\Biggr|\ge\delta\Sigma_n\Biggr\}
\le(\delta\Sigma_n)^{-2}
E \Biggl[ \Biggl\{ \sum^n_{i=1} \sigma_i (B_i-B_{i-1}) \Biggr\}^2
\Biggr]=\delta^{-2}.\quad
\end{equation}
Thus, by (\ref{eq:tnsip}), $\max_{1\le j\le n}|F_j|=\mathrm{O}_\mathrm{
p}(\Sigma_n+r_n\Delta_n)$. Observe that
%
\begin{equation}\label{eq:p1a}
\underline{V}^2_j-\Sigma^2_j =W_j-F_j^2/j, \qquad \mbox{where }   W_j
= \sum^j_{i=1} \sigma^2_i (e^2_i-1).
\end{equation}
By (\ref{eq:sip}), the same argument in
(\ref{eq:tnsip}) and (\ref{eq:maxin}) shows $W_j = \mathrm{O}_\mathrm{
p}(\Sigma^{*2}_n+r^*_n \Delta_n)$, uniformly. The desired
result then follows via (\ref{eq:p1a}).
\end{pf*}

\begin{pf*}{Proof of Theorem \protect\ref{thm:bootstrap}}
Denote by $\Phi(x)$ the standard normal distribution function.
By Proposition~\ref{cor:1} and Slutsky's theorem,
$\p(H_n/\hat\tau\le x)\to\Phi(x)$
for each fixed $x\in\R$. Since $\Phi(x)$ is a continuous
distribution, $\sup_{x\in\R}|\p(H_n/\hat\tau\le x)-\Phi(x)|=0$.
It remains to prove $\sup_{x\in\R}|\p^*(H^*_n\le x)-\Phi(x)|\to
0$, in probability. Notice that, conditioning on $\{e_i\}$,
$\{\xi_i\}$ are independent random variables with zero mean. By
the Berry--Ess\'een bound in Bentkus, Bloznelis and G\"{o}tze~\cite{bentkus}, there exists
a finite constant $c$ such that
%
\begin{equation}\label{eq:pfbta}
\sup_{x\in\R}|\p^*(H^*_n\le x)-\Phi(x)| \le c
\sum^n_{i=1} E^*(|\xi_i|^3) \Biggl\{\sum^n_{i=1} E^*(|\xi_i|^2) \Biggr\}^{-3/2},
\end{equation}
where $E^*$ denotes conditional expectations given $\{e_i\}$.
Clearly, $E^*(|\xi_i|^2)=\sigma^2_i e^2_i E(\alpha^2_1)$ and
$E(|\xi_i|^3)=\sigma^3_i |e^3_i| E(|\alpha^3_1|)$. Thus, under
the assumption $e_i\in\mathcal{L}^3$, we have $\sum^n_{i=1}
E^*(|\xi_i|^3)=\mathrm{O}_\mathrm{ p}(\sum^n_{i=1}\sigma^3_i)$. Meanwhile, by
the proof of Theorem~\ref{thm:0}, $\sum^n_{i=1}
E^*(|\xi_i|^2)=\sum^n_{i=1} \sigma^2_i e^2_i =
\{1+\mathrm{o}_\mathrm{ p}(1)\}\sum^n_{i=1}\sigma^2_i$. Therefore, the desired
result follows from (\ref{eq:pfbta}) in view of
(\ref{eq:bootcon}).
\end{pf*}

\begin{pf*}{Proof of Theorem \protect\ref{thm:test}}
For $cn\le j\le(1-c)n$, $c\le(1-j/n), j/n\le1-c$. For
$S_X(j)$ in (\ref{eq:SXj}), by (\ref{eq:thm0a}), we have
$\max_{cn\le j\le(1-c)n}|S_X(j)-\tau\widetilde{S}_X(j)|=\mathrm{O}_\mathrm{
a.s.}(r_n\Delta_n)$, where
\[
\widetilde{S}_X(j) = \biggl( 1-\frac{j}{n} \biggr)\sum^j_{i=1}\sigma_i
(B_i-B_{i-1}) -
\frac{j}{n} \sum^n_{i=j+1}\sigma_i (B_i-B_{i-1}).
\]
By (\ref{eq:thm0b}), $\max_{cn\le j\le(1-c)n}|(1-j/n)^2
\underline{V}^2_j + (j/n)^2 \overline{V}^2_j- V^2_j|=\mathrm{O}_\mathrm{
p}(\varpi_n)$, where
\[
V^2_j
= (1-j/n)^2 \sum^j_{i=1}\sigma^2_i + (j/n)^2 \sum^n_{i=j+1}\sigma^2_i
\quad \mbox{and}\quad  \varpi_n= (r^2_n\Delta^2_n + \Sigma^2_n)/n +
\Sigma^{*2}_n+r^*_n \Delta_n.
\]
For $cn\le j\le(1-c)n$, $V^2_j\ge c^2\Sigma^2_n$. Thus,
condition (\ref{eq:thmtestcon}) implies $\varpi_n=\mathrm{o}(V^2_j)$ and
$\{V^2_j+\mathrm{O}_\mathrm{ p}(\varpi_n)\}^{1/2}=V_j+\mathrm{O}_\mathrm{
p}(\varpi_n/V_j)$. Therefore, uniformly over $cn\le j\le
(1-c)n$,
\begin{eqnarray*}
T_n(j)-\tau\widetilde{T}_n(j)
=\frac{\tau\widetilde{S}_X(j)+\mathrm{O}_\mathrm{ a.s.}(r_n\Delta_n)}{V_j+\mathrm{O}_\mathrm{
p}(\varpi_n/V_j)}
- \frac{\tau\widetilde{S}_X(j)}{V_j} = \mathrm{O}_\mathrm{ p} \biggl\{\frac{r_n\Delta_n}{V_j}
+ \frac{\varpi_n \widetilde{S}_X(j)}{V_j^3} \biggr\}.
\end{eqnarray*}
By (\ref{eq:maxin}), $\max_{j} |\widetilde{S}_X(j)|=\mathrm{O}_\mathrm{ p}(\Sigma_n)$.
Thus, the result follows in view of $V_j\ge c \Sigma_n$.
\end{pf*}

\begin{pf*}{Proof of Theorem \protect\ref{thm:3}}
Condition $M_n\to0$ implies $\max_{1\le
j\le\ell_n}r(j)\Delta_n/\Sigma(j)\to0$. By (\ref{eq:thm0b}),
%
\begin{equation}\label{eq:omegaj}
\omega_j:=\frac{V^2(j)}{\Sigma^2(j)}-1
=\mathrm{O}_\mathrm{ p}\biggl\{\frac{\Sigma^{*2}(j) +r^*(j) \Delta_n}{\Sigma
^2(j)}+\frac{1}{k_n} \biggr\}
=\mathrm{O}_\mathrm{ p}(M_n)\to0.
\end{equation}
Define $U_j = \Sigma^{-1}(j) \sum_{i\in\mathcal{I}_j}\sigma_i
(B_i-B_{i-1})$. Clearly, $U_1,\ldots,U_{\ell_n}$ are
independent standard normal random variables. Thus, $\max_{1\le
j\le\ell_n}|U_j|=\mathrm{O}_\mathrm{ p}\{\sqrt{\log(\ell_n)}\} = \mathrm{O}_\mathrm{
p}\{\sqrt{\log(n)}\}$. By (\ref{eq:thm0a}),
$\underline{X}_n-\mu=\mathrm{O}_\mathrm{
p}\{(\Sigma_n+r_n\Delta_n)/n\}=\mathrm{O}_\mathrm{ p}(\Sigma_n/n)$. Recall
the definition of $D_j$ in (\ref{eq:Dj}). By the same argument
in (\ref{eq:thm0a}), using $\sqrt{1+x}=1+\mathrm{O}(x)$ as $x\to0$, we
have
\begin{eqnarray*}
D_j &=& \frac{k_n \{ \bar{X}(j) - \mu\}}{\Sigma(j)}\frac{1}{\sqrt
{1+\omega_j}}
+ \frac{k_n(\mu-\underline{X}_n)}{\Sigma(j)}\frac{1}{\sqrt{1+\omega_j}}
\\
&=& \biggl[ \tau U_j + \mathrm{O}_\mathrm{ a.s.} \biggl\{ \frac{r(j)\Delta_n}{\Sigma
(j)} \biggr\} \biggr] \{1+\mathrm{O}(\omega_j)\}
+\mathrm{O}_\mathrm{ p}\biggl\{\frac{k_n\Sigma_n}{n\Sigma(j)}
\biggr\}\\
&=& \tau U_j + \mathrm{O}_\mathrm{ p}\biggl\{ \sqrt{\log(n)} M_n + \frac{\Sigma
_n}{\ell_n \Sigma(j)} \biggr\}.
\end{eqnarray*}
By the latter expression and $\log(n)M_n\to0$, we can easily
verify $\hat\tau^2-\tau^2 =\mathrm{O}_\mathrm{ p}(\chi_n)$.
\end{pf*}
\end{appendix}

\section*{Acknowledgements}\vspace*{-2pt}
We are grateful to the associate editor and three anonymous referees
for their insightful comments that
have significantly improved this paper. We also thank Amanda Applegate
for help
on improving the presentation and Kyung-Ja Ha for
providing us the Seoul precipitation data. Zhao's research was
partially supported
by NIDA Grant P50-DA10075-15. The content is solely the
responsibility of the authors and does not necessarily
represent the official views of the NIDA or the NIH.


\printhistory

\end{document}